\documentclass[12pt]{amsart}

\usepackage{amsfonts, amssymb, amscd}
\usepackage{verbatim}
\usepackage{eucal}
\usepackage{amssymb}
\usepackage{mathrsfs}
\usepackage{graphicx}
\usepackage{psfrag}

\def\ii{{\rm i}}
\def\ee{{\rm e}}

\def\HH{{\mathcal H}}

\def\ZZ{{\mathbb Z}}
\def\RR{{\mathbb R}}

\def\CC{{\mathbb C}}
\def\QQ{{\mathbb Q}}

\def\mod{{\,{\rm mod}}}

\addtocounter{footnote}{1}

\newtheorem{lemma}{Lemma}[section]
\newtheorem{theorem}[lemma]{Theorem}
\newtheorem{corollary}[lemma]{Corollary}
\newtheorem{conjecture}[lemma]{Conjecture}
\newtheorem{proposition}[lemma]{Proposition}
\theoremstyle{definition}
\newtheorem{definition}[lemma]{Definition}

\newtheorem{remark}[lemma]{Remark}

\theoremstyle{remark}
\newtheorem*{proof*}{Proof}
\numberwithin{equation}{section}

\title{Holomorphic Eisenstein series with Jacobian twists}

\author{Lev A.  Borisov}
\address{Department of Mathematics \\ University of Wisconsin \\
  Madison \\ WI \\ 53706 \\ USA\\{\tt borisov@math.wisc.edu}}

\thanks{The author was partially supported by 
NSF grant DMS-0140172.}

\begin{document}

\begin{abstract}
For every point on the Jacobian of the modular 
curve $X_0(l)$ we define and study certain twisted
holomorphic Eisenstein series. These are particular cases
of a more general notion of twisted modular forms
which correspond to sections on the modular curve
$X_1(l)$ of the degree zero twists of line bundles of 
usual modular forms. We conjecture that a point
on the Jacobian is rational if and only if the ratios
of these twisted Eisenstein series of the same weights
have rational coefficients.
\end{abstract}

\maketitle

\section{Introduction}
Let $X_0(l)$ denote the modular curve of level $l$.
It is defined as the compactification of the quotient
of the upper half plane $\HH$ by the group $\Gamma_0(l)$
which consists of matrices $\Big(\begin{array}{cc}
a&b\\c&d\end{array}
\Big)\in SL_2(\ZZ)$ with $c\in l\ZZ$. 
We are interested in
the Jacobian of $X_0(l)$. We associate
to any point of this Jacobian a collection of holomorphic
functions $E_{i,k;h}(\tau)$ on the upper half plane $\HH$,
which
we call $h$-twisted Eisenstein series. Here $k\geq 3$
is an integer, $i\in \ZZ/l\ZZ$, and $h=h(\tau)$ is a modular form
of weight two for $\Gamma_0(l)$ which is related to
the point of the Jacobian by period integrals.

More specifically, to any weight two form $h$ 
we associate a unitary character of $\Gamma_0(l)$ by
\begin{equation}\label{char1}
\Big(\begin{array}{cc}
a&b\\c&d\end{array}
\Big)\mapsto {\rm exp}(2\pi\ii\Re\int_{\ii\infty}^{-\frac dc}h(s)\,ds)
\,.
\end{equation}
Then we use this unitary character to modify a definition 
of some Eisenstein series. The unitarity of the character
assures the convergence. It then turns out that the 
character and the series depend only on the point 
of the Jacobian, so one gets well-defined invariants 
of degree zero invertible sheaves on $X_0(l)$.

We remark that twisted Eisenstein series and 
modular forms of this paper should not be confused 
with twisted Eisenstein series of Goldfeld and Gunnells,
see \cite{Go} and \cite{GG}, which
involve only an additive version of the character \eqref{char1}.
However, a nonholomorphic analog of the Eisenstein
series considered in this paper has already appeared in
\cite{Petridis}.

The paper is organized as follows. In Section \ref{sec2}
we define the Eisenstein series 
twisted by modular forms $E_{k,i;h}(\tau)$ 
and prove their convergence. We use the Dirichlet summation 
trick to get a formula for the Fourier coefficients
of $E_{k,i;h}$. In Section \ref{sec3} we define a more general
class of functions, called $h$-twisted modular forms. The
definition is parallel to that of the usual modular forms, except
for the extra character \eqref{char1}. We give a geometric
interpretation of these twisted modular forms as sections
of the twists of the line bundles  of the usual modular forms  on 
the modular curves by degree zero line bundles, 
which justifies our terminology.
Section \ref{sec4} extends the Petersson inner product
to the twisted case. Section \ref{sec5} contains the
rationality Conjecture \ref{main}, as well as 
 some highly circumstantial
evidence for it. Finally, in Section \ref{sec6} we list
the open problems which currently by far  outnumber the results.

{\bf Acknowledgements.} I thank Paul
Gunnells and Gautam Chinta who read the first version of the 
paper and suggested several useful references. 

This paper is dedicated to the 60th birthday of 
Igor Dolgachev who was my PhD advisor at 
the University of Michigan in 1993-1996.  His love of
mathematics inspires everyone around him.
It is also a safe bet that he can pinpoint exactly
which 19th century mathematician has already
proved all the results of this paper.

\section{Eisenstein series twisted by modular forms}\label{sec2}
In this section we introduce the main objects of interest.
Let $l>0$ be a positive integer, which will be fixed throughout 
the paper. We will denote by $\tau$ the coordinate on the upper
half plane $\HH$ and will use the notation $q=\ee^{2\pi\ii\tau}$.
Let $h(\tau)=\sum_{n=1}^\infty a_n q^n$ 
be a cusp form of weight two for the congruence 
subgroup $\Gamma_0(l)$.
Let $H(\tau)=\int_{\ii\infty}^\tau  h(s)\,ds$ be the antiderivative of $h$. 

We first define the auxiliary series $\hat E_{k,\chi;h}(\tau)$.
\begin{definition}\label{hat}
For a positive integer $k$, a character $\chi:(\ZZ/l\ZZ)^*\to\CC$ 
and $h$ as above,
we define \emph{$h$-twisted Eisenstein series} by
$$
\hat E_{k,\chi;h}(\tau)=
\sum_{\Gamma_\infty\backslash \Gamma_0(l)}
(c\tau+d)^{-k} \chi(d) \ee^{2\pi\ii \Re H(-\frac dc)}.
$$
Here $\Gamma_\infty$ is the infinite cyclic subgroup generated by $\tau\to\tau+1$. 
\end{definition}

\begin{remark}\label{convention}
In Definition \ref{hat} and throughout the rest of the paper
we use the convention $H(-\frac d0)=0$. It is motivated
by the fact that $-\frac dc$ is the preimage of $\tau=\ii\infty$
under $\tau \to \frac {a\tau+b}{c\tau+d}$. So if $c=0$,
this preimage is $\ii\infty$, and $H(\ii\infty)=0$ by the 
definition.
\end{remark}

We can reformulate the above series as follows. The cosets of 
$\Gamma_0(l)$ over $\Gamma_\infty$ are given by pairs of integers
$(c,d)$ with $c\equiv 0\mod l$ and  ${\rm gcd}(c,d)=1$.
Hence, we have 
$$
\hat E_{k,\chi;h}(\tau)=
\sum_{c\in l\ZZ}\sum_{d\in\ZZ,{\rm gcd}(c,d)=1}
(c\tau+d)^{-k} \chi(d) \ee^{2\pi\ii \Re H(-\frac dc)}.
$$
We will be actually interested in a slight variation of this 
series, which is defined below.
\begin{definition}
We define 
$$
E_{k,i;h}(\tau)=
\sum_{c\in l\ZZ}\sum_{d\in\ZZ,{\rm gcd}(l,d)=1}
(c\tau+d)^{-k} \ee^{-2\pi\ii\frac {di}l} \ee^{2\pi\ii \Re H(-\frac dc)}.
$$
\end{definition}

\begin{proposition}\label{absolute}
For $k\geq 3$, the 
series $E_{k,i;h}(\tau)$ and $\hat E_{k,\chi;h}(\tau)$ absolutely converge to holomorphic
functions on the upper half plane.
\end{proposition}

\begin{proof}
The series of absolute values is 
$\sum_{(c,d)\in \ZZ^2-(0,0)} \vert c\tau+d\vert^{-k}$.
It has $\leq {(const)}R$ terms with $\vert c\tau+d\vert\in[R,R+1)$. 
Absolute convergence
on compacts is also clear, which implies that the resulting
sum is holomorphic.
\end{proof}

\begin{remark}
Proposition \ref{absolute} 
just barely fails for $k=2$ and fails quite
miserably for $k=1$, which is perhaps the case of most interest.
\end{remark}

By splitting the series $E_{k,i;h}$ according to ${\rm gcd}(c,d)$ we 
observe that it is a linear combination of series of type $\hat E$
with coefficients that 
involve the values of Dedekind L-functions at $k$ as well
as Gauss sums. To write it out explicitly, let us denote by $\chi_1,\ldots,\chi_{\phi(l)}$ the characters 
$(\ZZ/l\ZZ)^*\to \CC$.
\begin{proposition}\label{dedekind}
For $k\geq 3$ there holds 
$$
E_{k,i;h}(\tau) 
=\frac 1{\phi(l)}\sum_{j=1}^{\phi(l)}\Big(\sum_{t\in(\ZZ/l\ZZ)^*}
\ee^{-2\pi\ii\frac {ti}l}\chi_j^{-1}(t)\Big) L(\chi_j,k)\hat E_{k,\chi_j;h}(\tau)
$$
where $L(\chi,k)$ denotes the value of the 
Dedekind L-function at $k$, and $\phi$ is the Euler function.
\end{proposition}

\begin{proof}
For any $i\in \ZZ/l\ZZ$ and $d$ coprime to $l$, we have 
$$
\ee^{-2\pi\ii\frac {di}l} = \sum_{j=1}^{\phi(l)}
\chi_j(d) r_{i,j}
$$
where 
$$
r_{i,j}= \frac 1{\phi(l)} \sum_{t\in(\ZZ/l\ZZ)^*}
\ee^{-2\pi\ii\frac {ti}l}\chi_j^{-1}(t).
$$
This implies
$$
E_{k,i;h}(\tau) = 
\sum_{c\in l\ZZ}\sum_{d\in\ZZ,{\rm gcd}(l,d)=1}
\sum_{j=1}^{\phi(l)}
(c\tau+d)^{-k} 
r_{i,j}\chi_j(d) \ee^{2\pi\ii \Re H(-\frac dc)}
$$
$$
=\sum_{n\in\ZZ_{>0},gcd(n,l)=1}
\sum_{c_1\in l\ZZ}\sum_{d_1\in\ZZ,{\rm gcd}(c_1,d_1)=1}
\sum_{j=1}^{\phi(l)}
$$
$$
(c_1\tau+d_1)^{-k}n^{-k} 
r_{i,j}\chi_j(d_1)\chi_j(n) 
\ee^{2\pi\ii \Re H(-\frac {d_1}{c_1})}
$$
$$
=\sum_{c_1\in l\ZZ}\sum_{d_1\in\ZZ,{\rm gcd}(c_1,d_1)=1}
\sum_{j=1}^{\phi(l)}
(c_1\tau+d_1)^{-k}
r_{i,j}\chi_j(d_1)L(\chi_j,k)
\ee^{2\pi\ii \Re H(-\frac {d_1}{c_1})}
$$
since $L(\chi_j,k)=\sum_{n>0,gcd(n,l)=1}n^{-k}\chi_j(n)$.

So we have 
$$
E_{k,i;h}(\tau) = \sum_{j=1}^{\phi(l)} 
L(\chi_j,k) r_{i,j} \hat E_{k,\chi_j;h}(\tau)
$$$$=\frac 1{\phi(l)}\sum_{j=1}^{\phi(l)}\sum_{t\in(\ZZ/l\ZZ)^*}
L(\chi_j,k) \ee^{-2\pi\ii\frac {ti}l}\chi_j^{-1}(t) \hat E_{k,\chi_j;h}(\tau)
.
$$
The absolute convergence that allows us to change the order
of summation follows from the proof
of Proposition \ref{absolute}.
\end{proof}

We will now use the Dirichlet summation formula to 
find Fourier expansions of $E_{k,i;h}(\tau)$.
Namely, we have for $k\geq 2$ and $\Im x>0$
\begin{equation}\label{sum}
\sum_{n\in \ZZ} (x+n)^{-k}=
\sum_{m> 0} \frac{(-2\pi\ii)^{k}m^{k-1}}{(k-1)!}\ee^{2\pi\ii m x}.
\end{equation}

Equation \eqref{sum} allows us to rewrite the formula for $E_{k,i;h}(\tau)$.
\begin{proposition}\label{blah}
For $k\geq 3$ we have 
$$
E_{k,i;h}(\tau)
= \sum_{d\in\ZZ,gcd(d,l)=1}
d^{-k}\ee^{-2\pi\ii\frac {di}l}
+\frac{(-2\pi\ii)^{k}}{(k-1)!} \sum_{c\in l\ZZ_{>0}} 
\sum_{d_0\in (\ZZ/c\ZZ),
gcd(d_0,l)=1}$$
$$
c^{-k} 
\sum_{m> 0} m^{k-1} q^m\Big(
\ee^{-2\pi\ii\frac {d_0i}l}+(-1)^k\ee^{2\pi\ii\frac {d_0i}l}
\Big)
\ee^{2\pi\ii (m \frac {d_0}c+ \Re H(-\frac {d_0}c))}
$$
where  $q=\ee^{2\pi\ii\tau}$.
\end{proposition}

\begin{proof}
The terms with $c=0$ correspond to  the $q^0$ term.
Here we use the convention $H(-\frac d0)=0$ of Remark \ref{convention}.

For a given value of $c>0$, 
the possible values of $d$ are given
by $d=d_0+nc$ with $0<d_0<c,~gcd(d_0,l)=1$ and $n\in \ZZ$. 
For each $d_0$ the value of 
$\Re H(-\frac dc)$ is independent of $n$. 
Indeed, one easily sees that $H(\tau+1)=H(\tau)$, since $h(\tau)$
was a cusp form.
Then one uses \eqref{sum} to sum 
$(c\tau+d_0+nc)^{-k}=c^{-k}(\tau+\frac {d_0}c+n)^{-k}$ 
over all integer $n$.

Finally, we notice that terms for $c<0$ can be obtained for those
for $\vert c\vert >0$ by changing $d_0$ to $-d_0$.
Consequently, we sum for $c>0$ only 
by taking into account terms with $(-c)$.
\end{proof}

As a corollary, we observe that $E_{k,i;h}(\tau)$ has a
Fourier expansion.
\begin{corollary}\label{coeff}
For $k\geq 3$ one has
$$
E_{k,i;h}(\tau) = R_0 + \sum_{m>0}R_m q^m
$$
where 
$$
\begin{array}{rl}
R_m=& \frac{(-2\pi\ii)^{k}m^{k-1}}{(k-1)!} 
{\displaystyle\sum_{c\in l\ZZ_{>0}} \sum_{d_0\in (\ZZ/c\ZZ),
gcd(d_0,l)=1}}
\\
&
c^{-k} 
\Big(
\ee^{-2\pi\ii\frac {d_0i}l}+(-1)^k\ee^{2\pi\ii\frac {d_0i}l}
\Big)
\ee^{2\pi\ii (m \frac {d_0}c+ \Re H(-\frac {d_0}c))},~m>0\\
R_0=&\displaystyle\sum_{d\in\ZZ,gcd(d,l)=1}
d^{-k}\ee^{-2\pi\ii\frac {di}l}.
\end{array}
$$
\end{corollary}

\begin{proof}
Follows from  Proposition \ref{blah}.
\end{proof}

\begin{remark}
It is instructive to see what happens for $h=H=0$. 
Then for a fixed $c=c_0l$ the sum over $d_0$ that are the same $\mod l$
is zero unless $m$ is divisible by $c_0$. Indeed, for a
given $d_0\mod l$ we are summing over $d_0=s+rl\mod c_0l$
with $0\leq r<c_0$. The resulting sum of $\ee^{2\pi\ii m\frac r{c_0}}$
is zero unless $c_0\vert m$.
This gives
$$
R_m=\frac{(-2\pi\ii)^{k}m^{k-1}}{l^k(k-1)!} 
\sum_{c_0\vert m} c_0^{-k}
\sum_{\stackrel{d_0\in (\ZZ/c_0l\ZZ),}{gcd(d_0,l)=1}}
\Big(
\ee^{-2\pi\ii d_0(\frac il-\frac{m}{lc_0})}+(-1)^k
\ee^{2\pi\ii d_0(\frac il+\frac{m}{lc_0})}
\Big)
$$
$$
=
\frac{(-2\pi\ii)^{k}m^{k-1}}{l^k(k-1)!} 
\sum_{c_0\vert m} c_0^{-k}
c_0
\sum_{j\in (\ZZ/l\ZZ)^*}
\Big(
\ee^{-2\pi\ii j(\frac {i-\frac{m}{c_0}}{l})}+(-1)^k
\ee^{2\pi\ii j(\frac {i+\frac{m}{c_0}}l)}
\Big)
$$
$$
=\frac{(-2\pi\ii)^{k}}{l^k(k-1)!} 
\sum_{r\vert m} r^{k-1} 
\sum_{j\in (\ZZ/l\ZZ)^*}
\Big(
\ee^{-2\pi\ii \frac {j(i-r)}l}+(-1)^k
\ee^{2\pi\ii \frac {j(i+r)}l)}
\Big).
$$
This can be recognized as the coefficient by $q^m$ of 
a linear combination of the Eisenstein series from
\cite{highweight}.
Indeed, it is the sum over the divisors $r$ of $m$ of
an odd $(\mod l)$-polynomial function of $r$ of 
degree $(k-1)$, see \cite{highweight}.
In the particular case of prime $l$ we get
$$
\frac{(-2\pi\ii)^{k}}{l^{k-1}(k-1)!} 
\sum_{r\vert m} r^{k-1} 
\sum_{j\in (\ZZ/l\ZZ)^*}
\Big(
\delta_{r}^{i\mod l}+(-1)^k
\delta_{r}^{-i\mod l}-\frac 2l\delta_{k}^{0\mod 2}\Big)
$$
where $\delta$ is the Kronecker symbol.
\end{remark}

\begin{remark}\label{rational}
It is easy to see that the coefficients by $q^m$ 
of $\frac 1{(2\pi\ii)^k}E_{k,i;0}$ are 
rational for $m>0$. Indeed,
$\sum_{j\in (\ZZ/l\ZZ)^*}
\Big(
\ee^{-2\pi\ii \frac {j(i-r)}l}+(-1)^k
\ee^{2\pi\ii \frac {j(i+r)}l)}
\Big)
$ is invariant under the Galois group
of $\QQ[\ee^{2\pi\ii/l}]\supset \QQ$.
Since the space of Eisenstein series has a basis
whose elements have all Fourier coefficients in $\QQ$, the 
$q^0$ coefficients of $\frac {1}{(2\pi\ii)^k}E_{k,i;0}$ 
are rational as well. This is the main reason we prefer
$\frac {1}{(2\pi\ii)^k}E_{k,i;h}$ to their close relatives $
\frac {1}{(2\pi\ii)^k}\hat E_{k,\chi;h}$ that
do not have rational Fourier 
coefficients even in the untwisted case.
\end{remark}

Armed with the formulas for its Fourier coefficients, we can 
now try to define $E_{k,i;h}$ formally by their
Fourier expansions in the interesting 
cases $k=1$ and $k=2$. This approach works in the 
untwisted case, but fails in general. Indeed, let us analyze
the convergence of the series in Corollary \ref{coeff}.
We rewrite it as
$$
\sum_{c_0\in \ZZ_{>0}} \sum_{d_0\in (\ZZ/c_0l\ZZ),
gcd(d_0,l)=1}
c_0^{-k} 
\Big(
\ee^{-2\pi\ii\frac {d_0i}l}+(-1)^k\ee^{2\pi\ii\frac {d_0i}l}
\Big)
\ee^{2\pi\ii (m \frac {d_0}c+ \Re H(-\frac {d_0}{c_0l}))}.
$$
For a given $d_0\mod l=j\mod l$, we are  dealing with convergence
of 
$$
\sum_{c_0\in \ZZ_{>0}} \sum_{d_0\in (\ZZ/c_0l\ZZ),
d_0\equiv j \mod l}
c_0^{-k} 
\ee^{2\pi\ii (m \frac {d_0}c+ \Re H(-\frac {d_0}{c_0l}))}.
$$
There is no absolute convergence of the double sum,
even in the $k=2$ case. However, as the nontwisted case
$H=0$ 
suggests, one can try to look at the finite sums for fixed $c_0$
first and then look at the resulting series. It is conceivable 
that this series will converge  absolutely, at least for $k=2$.
Indeed, for small $h$ the values of the sums may be close 
to those for $h=0$ where all but a finite number are zero. 
But this is definitely not a proof. For instance, the values of 
$\Re H(-\frac {d_0}{c_0l})$ are not bounded, although 
one can show that they grow at most logarithmically in $c_0$.

\section{Twisted modular forms}\label{sec3}
In this section we put the definition of the $h$-twisted Eisenstein
series into a more general context of the $h$-twisted modular forms. We denote by $X_0(l)$ and $X_1(l)$ the modular
curves for the congruence subgroups 
$\Gamma_0(l)$ and $\Gamma_1(l)$. 
The group $\Gamma_1(l)$ is the subgroup of $\Gamma_0(l)$
whose diagonal elements are $1\mod l$.

\begin{definition}\label{t1}
An $h$-twisted modular form of weight $k$ with respect to 
$\Gamma_1(l)$ is a holomorphic
function $f$ on the upper half plane, 
which satisfies
\begin{enumerate}
\item For any $\Big(\begin{array}{cc}a&b\\c&d\end{array}\Big)
\in \Gamma_1(l)$ there holds
$$
f(\frac {a\tau+b}{c\tau+d}) = 
(c\tau+d)^k f(\tau) \ee^{-2\pi\ii \Re H(-\frac dc)}
.
$$
\item
For any $\Big(\begin{array}{cc}a&b\\c&d\end{array}\Big)
\in SL_2(\ZZ)$ the function
$$
f(\frac {a\tau+b}{c\tau+d})(c\tau+d)^{-k}
$$
is bounded near $\ii\infty$.
\end{enumerate}
\end{definition}

We immediately see that $E_{k,i;h}(\tau)$ and $\hat E_{k,\chi;h}(\tau)$
are $h$-twisted modular forms of level $l$. \begin{proposition}\label{Eismodular}
For $k\geq 3$, the functions $E_{k,i;h}(\tau)$ 
and $\hat E_{k,\chi;h}(\tau)$
are $h$-twisted modular form of weight $k$ for
$\Gamma_1(l)$.
\end{proposition}

\begin{proof}
We will first check the transformation properties of $E_{k,i;h}(\tau)$
under $\tau\to \frac {a\tau+b}{c\tau+d}$.
If the sum in the definition of $E_{k,i;h}$ is taken over 
$(c_1,d_1)$, then we can rewrite
the sum for $E_{k,i;h}(\frac{a\tau+b}{c\tau+d})$
in terms of the sum over $(c_2,d_2)=(c_1a+d_1c,c_1b+d_1d)$
as
$$
E_{k,i;h}(\frac{a\tau+b}{c\tau+d})=
\sum_{c_2\in l\ZZ}\sum_{\stackrel{d_2\in \ZZ,}{gcd(d_2,l)=1}}
(c_2\tau+d_2)^{-k} (c\tau+d)^k
\ee^{-2\pi\ii\frac {d_1i}l} \ee^{2\pi\ii \Re H(-\frac {d_1}{c_1})}.
$$
We notice that $d_1\equiv d_2\mod l$, so the first exponential term
equals $\ee^{-2\pi\ii\frac {d_2i}l}$. It then suffices to show that
$H(-\frac {d_1}{c_1})=H(-\frac {d_2}{c_2})-H(-\frac dc)$.
Since $H(\tau)$ is the antiderivative of a weight two $\Gamma_0(l)$
modular form $h$, the transformation properties of $h(\tau)$ imply
$$
H(\frac {-d\tau+b}{c\tau-a})=H(\tau)+ C
$$
where $C$ is independent of $\tau$.  By plugging in $\tau=\ii\infty$,
we calculate $C$ to get
$H(\frac {d_1\tau+b_1}{c_1\tau-a_1})=H(\tau)+ H(-\frac {d}{c})$.
Then we plug in $\tau = -\frac {d_1}{c_1}$ to get 
$H(-\frac {d_2}{c_2})=H(-\frac {d_1}{c_1})+H(-\frac dc)$.

We now check that for any 
$\Big(\begin{array}{cc}a&b\\c&d\end{array}\Big)
\in SL_2(\ZZ)$ the function
$$
E_{i,k;h}(\frac {a\tau+b}{c\tau+d})(c\tau+d)^{-k}
$$
is bounded near $\ii\infty$.
We have
$$
E_{k,i;h}(\frac{a\tau+b}{c\tau+d})(c\tau+d)^{-k}=
\sum_{c_1\in l\ZZ}\sum_{\stackrel{d_1\in \ZZ,}{gcd(d_1,l)=1}}
(c_2\tau+d_2)^{-k} 
\ee^{-2\pi\ii\frac {d_1i}l} \ee^{2\pi\ii \Re H(-\frac {d_1}{c_1})}.
$$
Each term is bounded by $\vert c_2\tau+d_2\vert^{-k}$.
We extend the summation set to sum over $(c_2,d_2)\in \ZZ^2-\{(0,0)\}$
and assume $\Im\tau\geq 1$. Then $\vert c_2\tau+d_2\vert\geq 1$,
and we consider the subsums for 
$\vert c_2\tau+d_2\vert \in [m,m+1)$ for all positive integer $m$.
If we slice the annulus $\vert z\vert \in [m,m+1)$ into $m$ equal
pieces according to the polar angle, we see that each slice contains
at most a constant number of elements of $\ZZ\tau+\ZZ$, since 
this lattice has no elements of length less than one. 
Importantly, this constant $C$ is independent of $\tau$.
This gives
$$
\sum_{(c_2,d_2)}\vert c_2\tau+d_2\vert^{-k}\leq C\sum_{m>0}m^{-k+1},
$$
which converges for $k\geq 3$.

The argument for $\hat E_{k,\chi;h}$ is similar and is left
to the reader.
\end{proof}

We can extend the Proposition \ref{Eismodular}
to describe the transformation
properties of $E_{k,i;h}$ and $\hat E_{k,\chi;h}$
under $\Gamma_0(l)$.
\begin{proposition}
For $k\geq 3$, 
for every $\tau\to\frac {a\tau+b}{c\tau+d}\in \Gamma_0(l)$
one has
$$
E_{k,i;h}(\frac{a\tau+b}{c\tau+d})=E_{k,ai;h}(\tau)(c\tau+d)^k
   \ee^{-2\pi\ii \Re H(-\frac dc)}
$$
$$
E_{k,\chi;h}(\frac{a\tau+b}{c\tau+d})=E_{k,\chi;h}(\tau)(c\tau+d)^k
   \chi(a)\ee^{-2\pi\ii \Re H(-\frac dc)}.
$$
\end{proposition}

\begin{proof}
The argument for $E_{k,i;h}$ follows that of Proposition \ref{Eismodular}.
The only difference
is what happens to the term $\ee^{-2\pi\ii \frac {d_1i}l}$, where
one now observes that $d_2\equiv d_1d\mod l$, so $d_1\equiv ad_2\mod l$.

The argument for $\hat E_{k,\chi;h}$ is similar and is left to
the reader.
\end{proof}

\begin{remark}\label{neben}
The space of $h$-twisted modular forms of weight $k$ and 
with respect to $\Gamma_1(l)$ admits a natural
action of $\Gamma_0(l)/\Gamma_1(l)\cong (\ZZ/l\ZZ)^*$. It sends
$f(\tau)$ to
$$
f(\frac{a\tau+b}{c\tau+d})(c\tau+d)^{-k}\ee^{2\pi\ii \Re H(-\frac dc)}.
$$
The eigenvectors of 
this action will be called $h$-twisted modular forms of nebentypus $\chi$.
They satisfy 
$$
f(\frac{a\tau+b}{c\tau+d})=f(\tau)(c\tau+d)^k \chi(a)
\ee^{-2\pi\ii \Re H(-\frac dc)}
$$
for $\Big(\begin{array}{cc}a&b\\c&d\end{array}\Big)$ in
$\Gamma_0(l)$.
We notice that $\hat E_{k,\chi;h}$ 
have nebentypus $\chi$.
\end{remark}

\begin{remark}
We will call the nebentypus $\chi=1$ the $h$-twisted 
$\Gamma_0(l)$ modular forms. For example, $E_{k,0;h}(\tau)$
is one such form for $k\geq 3$. As in the usual case, the action
of $-{\bf 1}\in \Gamma_0(l)$ causes the spaces of odd weight
to be zero.
\end{remark}

It is clear that 
for a given $h$, the space of 
$h$-twisted modular forms for $\Gamma_1(l)$ is a graded module over
the ring of modular forms for $\Gamma_1(l)$. In what follows we give
a geometric interpretation of this space. We can safely assume $l\geq 3$
so the action of $\Gamma_1(l)$ on the upper half plane is free.
We denote by $X_1(l)$ the corresponding modular curve which is obtained
by the compactification of its open subset $\HH/\Gamma_1(l)$.
There is a line bundle $\mathcal L$ of weight one modular forms on 
$X_1(l)$, such that the space $M_k(\Gamma_1(l))$ 
of usual modular forms of weight $k$ is naturally identified with
the space of global sections $H^0(X_1(l),{\mathcal L}^{\otimes k})$.

\begin{proposition}
For each $h$ there is a degree zero 
invertible ${\mathcal F}={\mathcal F}(h)$
such that the space of weight $k$ $h$-twisted modular forms of level $n$
is naturally identified with $H^0(X_1(l),{\mathcal F}\otimes
{\mathcal L}^{\otimes k})$.
\end{proposition}

\begin{proof}
For an $h$-twisted modular form $f$ we
define the order of zero of $f$ at a cusp to be 
the exponent in the corresponding Fourier expansion.
Combined with the zeroes in the open set $\HH/\Gamma_1(l)$,
they give a divisor on $X_1(l)$ which we denote by
$Div(f)$. 

We observe that the divisor class of $Div(f)$ 
depends only on $h$ and $k$. Indeed, the ratio 
$f(\tau)=\frac {f_1(\tau)}{f_2(\tau)}$ of two weight $k$ 
$h$-twisted modular forms is a $\Gamma_1(l)$ modular function.
It remains to observe that $Div(f)=Div(f_1)-Div(f_2)$.
More generally, for $h$-twisted forms $f_{i}$ of weights $k_i$ 
one has $Div(f_1)-Div(f_2)\sim (k_1-k_2) Div(L)$ where $Div({\mathcal L})$
is the divisor of any (meromorphic) modular form of weight $1$.
We define $\mathcal F$ as the invertible sheaf
that corresponds to $Div(E_{4,0;h})-Div(E_{4,0;0})$.

Every $h$-twisted modular form of weight $k$ gives 
a global section of ${\mathcal F}\otimes {\mathcal L}^{\otimes k}$.
Conversely, given a global section of 
${\mathcal F}\otimes {\mathcal L}^{\otimes k}$,
one can write it as a product of  $E_{4,0;h}(\tau)$
and a modular function on $X_1(l)$. This product is easily seen
to be a holomorphic function on $\HH$ that satisfies
the definition of the $h$-twisted modular form of weight $k$.

It remains to see that the degree of $\mathcal F$ is zero. 
Consider a ratio of an $h$-twisted and a usual form of the same weight
$k$, which we can assume to be nonzero at the cusps.
It is a meromorphic function $g(\tau)$ on the upper half plane, which is 
bounded at the cusps.
To find the number of its zeroes and poles in the fundamental domain,
one needs to integrate to find the sum of residues of 
 $d\log g(\tau)$ in it.
Transformation properties of $g$ imply that 
$d\log g(\tau)$ is a meromorphic differential form on $X_1(l)$.
Consequently, the sum of its residues is zero, which leads to
 $\deg {\mathcal F}=0$, as desired.
\end{proof}

\begin{proposition}
In the notations above, the invertible sheaf ${\mathcal F}(h)$
is invariant under the action of $\Gamma_0(l)/\Gamma_1(l)$. 
Moreover, it admits a natural linearization for this action.
\end{proposition}

\begin{proof}
The key point here is that $h(\tau)$ was a 
$\Gamma_0(l)$ modular form.
The linearization is provided by the action of 
$\Gamma_0(l)/\Gamma_1(l)$ from Remark \ref{neben}.
\end{proof}

\begin{remark}\label{F0}
We may as well talk about the invertible sheaf ${\mathcal F}_0=
{\mathcal F}_0(h)$
of degree zero on the modular curve $X_0(l)$. 
 It can be either defined as the 
$\Gamma_0(l)/\Gamma_1(l)$-invariant part of the pushforward
$\pi_*\mathcal F$ for the quotient map $\pi:X_1(l)\to X_0(l)$, 
or as a sheaf that corresponds to
the divisor $Div(E_{4,0;h})-Div(E_{4,0;0})$ on $X_0(l)$.
We have $\pi^*{\mathcal F}_0\cong {\mathcal F}$.
\end{remark}

\begin{remark}
The operation of tensoring by an invertible sheaf is commonly
referred to as twisting, which justifies our terminology.
Unfortunately, this term is already used in \cite{GG} in
the context of Eisenstein series. We hope that 
this does not lead to a confusion, since the series we construct
in this paper are quite different from those in \cite{GG}.
\end{remark}

\begin{theorem}\label{surj}
Any degree zero invertible sheaf
${\mathcal F}_0$ on $X_0(l)$ 
is isomorphic to ${\mathcal F}_0(h)$ of Remark \ref{F0}
for some weight two 
$\Gamma_0(l)$ cusp form $h$.
\end{theorem}

\begin{proof}
The map $h\mapsto {\mathcal F}_0(h)$ 
can be described  in terms of the periods.
Let $r$ be the genus of $X_0(l)$, which we assume to 
be positive. Denote by $K_{X_0(l)}$ the
sheaf of holomorphic $1$-forms on $X_0(l)$.
The Jacobian  of $X_0(l)$ is isomorphic to the
quotient of $\CC^{r}=H^0(X_0(l),K_{X_0(l)})^\vee$ 
by the lattice of periods, defined as the image
of $H_1(X_0(l),\ZZ)$ via integration, see \cite{GH}.
We would like to calculate the point of the Jacobian
that corresponds to ${\mathcal F}_0(h)$ in these
terms. This means that for any 
$g\in H^0(X_0(l),K_{X_0(l)})$ we would like to find
the integral of $g$ from zeroes of some $h$-twisted 
weight $k$ form $f_{h}$
to zeroes of some usual 
weight $k$ form $f$, up to the lattice of periods.

By an easy Riemann-Roch calculation, we can safely
assume that $f_h$ and $f$ do not 
vanish at the extra symmetry points or cusps,
where by  extra symmetry points we mean those 
that are equivalent to $\ii$ or $\ee^{2\pi\ii/3}$ modulo
$SL_2(\ZZ)$.
Let $G(\tau)=\int_{\ii\infty}^{\tau}g(s)\,ds$ be an
antiderivative of $g$. In order to find the integral of 
$g$ from the zeroes of $f$ to the zeroes of $f_h$,
we need to find the integral
$$
I=\frac 1{2\pi\ii}\int_{\partial \HH/\Gamma_0(l)}
G(\tau)
\Big(
\frac {f'_h(\tau)}{f_h(\tau)}
-
\frac {f'(\tau)}{f(\tau)}
\Big)
\,d\tau
$$
over the boundary of a fundamental domain of
$\Gamma_0(l)$ action on $\HH$. Different choices
of fundamental domains will give answers that differ
by periods of $g$. 

\begin{figure}[tbh]
\begin{center}
\includegraphics[scale = .3]{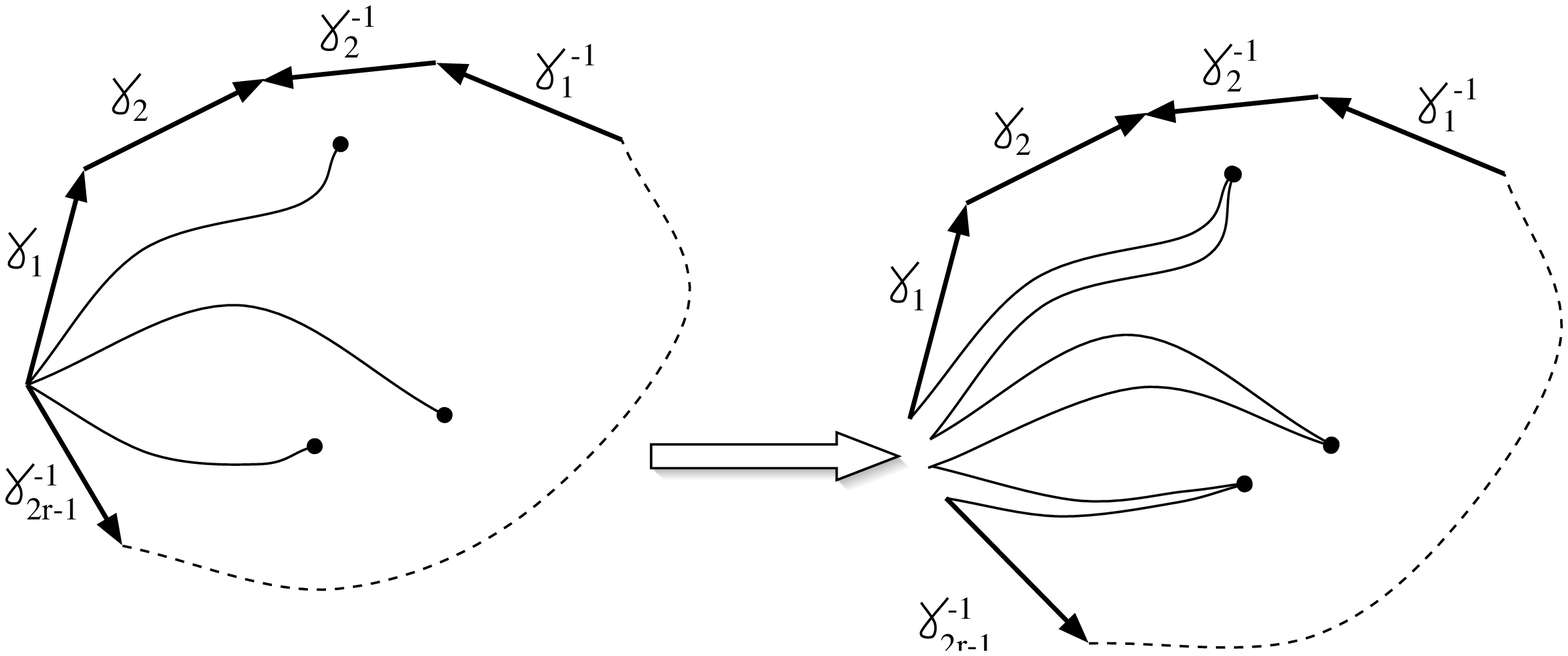}
\end{center}
\caption{\label{fig1}}
\end{figure}

Let us pick a fundamental domain as follows.
Pick a point $p\in X_0(l)$ and cut $X_0(l)$ 
to represent it by a $4r$-gon with pairs of sides identified
in the usual manner, such that all vertices map to $p$.
We will assume that the cusps, 
the points of extra symmetry, as well as the zeroes of 
$f_h$ and $f$ 
lie in the interior of this polygon.
We will draw cuts from a vertex of 
the $4r$-gon to the cusps and extra symmetry points,
see Figure \ref{fig1}.
We can represent the resulting $4r$-gon without these
extra cuts by a region in the 
upper half plane by picking a preimage of  a point in
the middle and extending all the paths from it by continuity.
This will have an effect of opening up the cuts to the cusps
and to those extra symmetry points whose isotropy subgroup
in $\Gamma_0(l)$ is bigger than just $\{\pm 1\}$,
see Figure \ref{fig1}.

The boundary arcs of the resulting fundamental domain 
$\mathcal D$ are paired up according
to the gluing required to make $X_0(l)$. These arcs
will then be related by a coordinate transformation
$\tau\to \frac {a\tau+b}{c\tau+d}$ from $\Gamma_0(l)$.
We first investigate pairs of paths that go
to the cusps or symmetry points of $X_0(l)$. Whether or
not these cuts open up, the values 
of $G(\tau)
\Big(
\frac {f'_h(\tau)}{f_h(\tau)}
-
\frac {f'(\tau)}{f(\tau)}
\Big)
$
at the opposite sides of the cut are the same.
Indeed, $G(\frac {a\tau+b}{c\tau+d})=G(\tau)$ 
for any $\tau\to \frac {a\tau+b}{c\tau+d}$ from $\Gamma_0(l)$
that has a fixed point. Similar statement holds for 
$H(\tau)$, and the transformation properties of 
$f_h$ and $f$ show that $d\log f_h- d\log f$ is
unchanged under such transformations.

So we are left to integrate over the arcs $p_j\to q_j$
of the boundary that correspond to the sides of the original
$4r$-gon. Pair them up according to the gluings.
The two arcs $\gamma_j$ and $\gamma_j^{-1}$ (see Figure 
\ref{fig1})
differ by $\tau\to \frac {a_j\tau+b_j}{c_j\tau+d_j}$
from $\Gamma_0(l)$, and we are looking at
$$
I=\frac 1{2\pi\ii}\sum_{j=1}^{2r} \int_{p_j}^{q_j}
\Big(
G(\frac {a_j\tau+b_j}{c_j\tau+d_j})
\Big(
\frac {f'_h(\frac {a_j\tau+b_j}{c_j\tau+d_j})}
{f_h(\frac {a_j\tau+b_j}{c_j\tau+d_j})}
-\frac {f'(\frac {a_j\tau+b_j}{c_j\tau+d_j})}
{f(\frac {a_j\tau+b_j}{c_j\tau+d_j})}
\Big)(c_j\tau+d_j)^{-2}
$$$$
-
G(\tau)
\Big(
\frac {f'_h(\tau)}
{f_h(\tau)}
-\frac {f'(\tau)}
{f(\tau)}
\Big)
\Big)
\,d\tau
$$
Transformation properties of $G$, $f_h$ and $f$ 
allow us to simplify the integrand to 
be 
$$
\Big(
(G(\tau) + \alpha_j)
\Big(\frac {f'_h(\tau)}
{f_h(\tau)}
-\frac {f'(\tau)}
{f(\tau)}
\Big)
-
G(\tau)
\Big(
\frac {f'_h(\tau)}
{f_h(\tau)}
-\frac {f'(\tau)}
{f(\tau)}
\Big)
\Big)
$$
$$
=\alpha_j
\Big(
\frac {f'_h(\tau)}
{f_h(\tau)}
-\frac {f'(\tau)}
{f(\tau)}
\Big)
$$
where $\alpha_j=G(\frac{a_j\tau+b_j}{c_j\tau+d_j})-G(\tau)$.
Hence we get 
$$
I=\frac 1{2\pi\ii}\sum_{j=1}^{2r}
\int_{p_j}^{q_j}\alpha_j(\Big(
\frac {f'_h(\tau)}
{f_h(\tau)}
-\frac {f'(\tau)}
{f(\tau)}
\Big)\,d\tau
=\frac 1{2\pi\ii}\sum_{j=1}^{2r}
\alpha_j
\log\Big( \frac {f_h(\tau)}{f(\tau)}\Big)
\Big\vert_{p_j}^{q_j}
$$$$
=
-\sum_{j=1}^{2r} \alpha_j\Re H(-\frac {\hat d_j}{\hat c_j})\mod ~{periods},
$$
where $\tau\to 
\frac {\hat a_j\tau+\hat b_j}{\hat c_j \tau+\hat d_j}$
is an element of $\Gamma_0(l)$
that sends $p_j\mapsto q_j$. 
Indeed, the transformation properties of $f_h$ and $f$
lead to the transformation properties of the logarithms, up
to integer constants that are independent of $g$, and $\alpha_i$
are the period integrals.

We observe that $\alpha_j$ are the period integrals
of $g(\tau)\,d\tau$ over the standard symplectic basis of 
$H_1(X_0(l),\ZZ)$, as are $H(-\frac {\hat d_j}{\hat c_j})$
are the integrals of $h(\tau)\,d\tau$ over the same
basis. Indeed, 
$H(\frac {\hat a_j\tau+\hat b_j}{\hat c_j \tau+\hat d_j})
=H(\tau) - H(-\frac {\hat d_j}{\hat c_j})$ in view of $H(\ii\infty)=0$,
so $-H(-\frac {\hat d_j}{\hat c_j})  = H(q_j)-H(p_j)=
\int_{p_j}^{q_j}h(\tau)\,d\tau$. To show that 
any point of the Jacobian can be given by ${\mathcal F}_h$,
it is enough to show that the pairing (of real vector spaces)
$$
H_1(X_0(l),\RR) \times H^0(X_0(l),K_{X_0(l)})\to \RR,~
(\gamma, w)\mapsto \Re \int_\gamma w
$$
is nondegenerate. This is a standard general fact, which 
is true for any projective curve $X$, but we provide
the argument below for the benefit of the reader.

By picking a basis $(w_1,\ldots,w_r)$ dual to
a half of a symplectic basis $(\gamma_1,\ldots,\gamma_{2r})$, 
the period matrix $\Omega=(\int_{\gamma_j}w_i)$
is given by $({\bf I}_r,Z)$, where $\Im Z$ is positive definite
(see \cite{GH}). If $\int_{\gamma_j} w$ is purely imaginary
for all $j$, then we first see that $w=\sum_{i=1}^r \lambda_i w_i$
with purely imaginary $\lambda_i$. Then $(\lambda_i)$
lies in the kernel of $Z$, so $\Im Z>0$ implies 
$(\lambda_i)={\bf 0}$.
\end{proof}

\begin{corollary}\label{func}
Two weight two cusp forms $h_1$ and $h_2$ 
have ${\mathcal F}_1\cong {\mathcal F}_2$ if
and only if  $\Re(H_1(-\frac dc)-H_2(-\frac dc))\in\ZZ$ 
for all $c\in l\ZZ$ and $d\in\ZZ,gcd(d,l)=1$.
In particular, the Eisenstein series $E_{k,i;h}$ depend
only  on the point of the Jacobian.
\end{corollary}

\begin{proof}
The \emph{if} part is clear, since in this case the Eisenstein
series are identical.

To see the \emph{only if} part, we observe that in the proof
of Theorem \ref{surj}, the sheaf $\mathcal F$ is given by 
a linear combination of the periods with coefficients
$\Re H(-\frac {\hat d_j}{\hat c_j})$, $j=1,\ldots, 2r$. If 
${\mathcal F}_1\cong {\mathcal F}_2$, then 
\begin{equation}\label{bas}
\Re(H_1(-\frac {\hat d_j}{\hat c_j})-H_2(-\frac {\hat d_j}{\hat c_j}))\in\ZZ
\end{equation}
for $j=1,\ldots,2r$. For any $\Big(\begin{array}{cc}
a&b\\c&d\end{array}\Big)$ in $\Gamma_0(l)$, the
values of $\Re H_1(-\frac dc)- \Re 
H_2(-\frac dc)$ are integer linear combinations of 
the periods in \eqref{bas}. Indeed, this follows from
the fact that \eqref{bas} gives integrals of 
$(h_1(\tau)-h_2(\tau))\,d\tau$
over a basis of $H_1(X_0(l),\ZZ)$.
Finally, any ratio $-\frac dc$ with $c\in l\ZZ$ and 
$gcd(d,l)=1$ can be obtained from a $\Gamma_0(l)$
matrix, after one cancels off common factors of $c$ and $d$.
\end{proof}

\begin{remark}
Corollary \ref{func} shows that 
each Fourier coefficient of $E_{k,i;h}$ gives
a function on the Jacobian of $X_0(l)$. However,
these functions are clearly not holomorphic in general,
else they would have to be constant.
\end{remark}

\begin{remark}
An easy application of the Riemann-Roch theorem 
shows that the dimension of the space of $h$-twisted
modular forms of weight $k\geq 2$ coincides with that 
of the space of usual modular forms of weight $k$.
This holds for either $\Gamma_0(l)$ or $\Gamma_1(l)$,
as well as for any nebentypus.
\end{remark}

\section{Petersson inner product for twisted modular 
forms and twisted Eisenstein series}\label{sec4}

In this section we define the Petersson inner product
on the space of $h$-twisted modular forms of given weight $k$
and investigate its basic properties.

\begin{definition}
A $\Gamma_1(l)$-modular $h$-twisted form is called a \emph{cusp form},
if in the second condition of Definition \ref{t1} the function
$$
f(\frac {a\tau+b}{c\tau+d})(c\tau+d)^{-k}
$$
has limit $0$ as $\tau\to\ii\infty$ for all 
$\Big(\begin{array}{cc}a&b\\c&d\end{array}\Big)$.
\end{definition}

\begin{definition}\label{Petersson}
Let $f$ and $g$ be two 
$h$-twisted modular forms for $\Gamma_1(l)$
of weight $k$. Assume that at least one of then is
a cusp form. 
We define Petersson inner product $\langle f,g\rangle$
by
$$
\langle f,g\rangle=\int_{\Gamma_1(l)\backslash{\mathcal H}}
f(\tau)\overline{g(\tau)} y^{k-2}\,dxdy
$$
where $\tau=x+\ii y$. The convergence is assured by 
the cusp condition.
\end{definition}

\begin{remark}
We have to check that the above definition is independent of 
the choice of the fundamental domain. If
$\tau = \frac {a\tau_1+b}{c\tau_1+d}$, then
$$2\ii y= \frac {a\tau_1+b}{c\tau_1+d}-
\frac {a\bar \tau_1+b}{c\bar\tau_1+d}
=\frac {(a\tau_1+b)(c\bar\tau_1+d)-(c\tau_1+d)(a\bar \tau_1+b)}
{(c\tau_1+d)(c\bar\tau_1+d)},
$$
$$
\frac {\tau_1-\bar \tau_1}{(c\tau_1+d)(c\bar\tau_1+d)}=
\frac {2\ii y_1}{(c\tau_1+d)(c\bar\tau_1+d)}
$$
and
$$
dxdy = \frac 1{(c\tau_1+d)^2(c\bar\tau_1+d)^2} \,dx_1dy_1
$$
so
$$
f(\tau)\overline{g(\tau)} y^{k-2}\,dxdy=
f(\frac {a\tau_1+b}{c\tau_1+d})\overline{g(\frac {a\tau_1+b}{c\tau_1+d})}
\frac {y_1^{k-2}}{(c\tau_1+d)^{k}(c\bar\tau_1+d)^{k}}
\,dx_1dy_1
$$
$$
=f(\tau_1)
\ee^{-2\pi\ii \Re H(-\frac dc)}\overline{g(\tau_1)}
\ee^{2\pi\ii {\Re H(-\frac dc)}}{y_1^{k-2}}
\,dx_1dy_1=f(\tau_1)\overline{g(\tau_1)}
{y_1^{k-2}}
\,dx_1dy_1.
$$
\end{remark}

\begin{remark}
It is clear that the above defined Petersson pairing 
restricts to a nondegenerate Hermitean form on
the space of $h$-twisted cusp forms. 
\end{remark}

The following proposition extends the orthogonality of
Eisenstein series and cusp forms (see \cite{Lang})
to the $h$-twisted case,
for the particular class of Eisenstein series considered
in this paper.
\begin{proposition}\label{orthogonal}
For $k\geq 3$ we have
$$\langle E_{k,i;h},g\rangle=0$$
for any weight $k$ $h$-twisted cusp form $g$ for $\Gamma_1(l)$.
\end{proposition}

\begin{proof}
It is enough to check that $\langle\hat E_{k,\chi;h},g\rangle =0$
for all characters $\chi:(\ZZ/l\ZZ)^*\to \CC^*$.
We can also assume that $g$ is an eigenform for $\Gamma_0$-action, as in Remark \ref{neben}. Denote the character of 
$g$ by $\psi$.
We have 
$$
\langle \hat E_{k,\chi;h},g\rangle=
\int_{\Gamma_1(l)\backslash{\mathcal H}}
\sum_{\Gamma_\infty\backslash \Gamma_0(l)}
(c\tau+d)^{-k}\chi(d)
\ee^{2\pi\ii\Re H(-\frac dc)}\overline{g(\tau)} y^{k-2}\,dxdy.
$$
We switch the order of integration and summation
and then switch to a different domain for each term.
Namely, we make a substitution $\tau=\frac {d\tau_1-b}{-c\tau_1+a}$
for each term of the $\hat E$-series. This gives
$$
\langle \hat E_{k,\chi;h},g\rangle=
\sum_{\Gamma_\infty\backslash \Gamma_0(l)}
\int_{\Gamma_1(l)\backslash{\mathcal H}}
(-c\tau_1+a)^k 
\chi(d)
\ee^{2\pi\ii\Re H(-\frac dc)} 
\cdot$$$$\cdot
(-c\bar \tau_1+a)^k \ee^{-2\pi\ii \Re H(-\frac dc)}
\overline{\psi(d)}\,\overline{g(\tau_1)}
\frac {y_1^{k-2}dx_1dy_1}
{(-c\tau_1+a)^k(-c\bar\tau_1+a)^k}
$$$$=\sum_{\Gamma_\infty\backslash \Gamma_0(l)}
\int_{\Gamma_1(l)\backslash{\mathcal H}}
\chi(d)\overline{\psi(d)}\,
\overline{g(\tau_1)}
y_1^{k-2}\,dx_1dy_1.
$$
This allows us to rewrite the pairing as
$$
\sum_{\Gamma_1(l)\backslash\Gamma_0(l)}
\chi(d)\overline{\psi(d)}
\int_{\Gamma_\infty\backslash\mathcal H}
\overline{g(\tau_1)}
y_1^{k-2}\,dx_1dy_1.
$$
We then observe that  $\Gamma_\infty\backslash{\mathcal H}$
can be thought of as $0\leq x\leq 1,y\geq 0$ and for a fixed
$y$
$$
\int_{\ii y}^{1+\ii y} \overline{g(\tau)}\,dx
=\sum_{m>0} \bar r_m\ee^{-2\pi\ii my}\int_{x=0}^{x=1}
\ee^{-2\pi\ii m x}\,dx = 0.
$$
At the last step we used that $g$ has a Fourier expansion
with a zero leading term. Indeed, we have $H(\ii\infty)=0$,
which implies that $g(\tau+1)=g(\tau)$, and the $g(\ii\infty)=0$
follows from $g$ being a cusp form.
\end{proof}

\section{Rationality conjecture}\label{sec5}
The motivation behind our investigation of 
twisted Eisenstein series is that 
they might provide a tool for studying the Jacobian of
the modular curve $X_0(l)$. We now recall that $X_0(l)$
comes equipped with the natural $\QQ$-structure,
which then gives a $\QQ$-structure on its Jacobian.
The following conjecture relates the rationality of 
the coefficients of  twisted Eisensten series and the 
rationality of the point on the Jacobian.
\begin{conjecture}\label{main}
The sheaf ${\mathcal F}_0(h)$ on $X_0(l)$ gives
a rational point on the Jacobian if and only if 
for each  $k$ there exists a nonzero formal power
series $R_k(q)$ such that all  coefficients 
of the $q$-series $E_{k,i;h}R_k(q)$ are \emph{rational}
for all $i$. Equivalently, the ratios of $E_{k,i;h}$ and
$E_{k,j;h}$ have rational Fourier coefficients for
all $k$, $i$ and $j$, provided $E_{k,j;h}\not\equiv 0$.
\end{conjecture}

\begin{remark}
A stronger version of Conjecture \ref{main} would
require $R_k=R^k$ for some series $R=R(q)$.
\end{remark}

We devote the remainder of the section to the evidence in favor of Conjecture \ref{main}.
Since the evidence is mostly circumstantial anyway,
we will be somewhat sketchy and leave the details to
the reader.

First of all, the \emph{if} part of the conjecture holds 
up to a finite index subgroup as 
long as $E_{k,i;h}$ do not have common zeroes for
some $k>0$. Indeed, the existence of $R_k$ implies 
that the module over the ring of $\Gamma_1(l)$ modular
forms which is generated by $E_{k,i;h}$ is defined 
over $\QQ$. By the usual correspondence between graded
modules and sheaves we see that the invertible sheaf
${\mathcal F}(h)\otimes {\mathcal L}^{\otimes k}$ is
defined over $\QQ$. Since $\mathcal L$ is defined over $\QQ$,
we see that the invertible sheaf  ${\mathcal F}(h)$ 
on $X_1(l)$ is defined over $\QQ$. 
Denote $\pi:X_1(l)\to X_0(l)$ and recall that 
${\mathcal F}(h)=\pi^*{\mathcal F}_0(h)$.
Then 
$$
{\mathcal F}_0^{\phi(l)/2}\cong (\Lambda^{\phi(l)/2}
\pi_*{\mathcal F}(h))\otimes (\Lambda^{\phi(l)/2}\pi_*
\pi^*{\mathcal O})^{-1}
$$
is defined over $\QQ$.

For the \emph{only if} part of the conjecture, observe 
that it holds for $h=0$, in view of Remark \ref{rational}.
In general, for every rational point on the Jacobian 
of $X_1(l)$ that 
corresponds to a sheaf ${\mathcal F}=\pi^*{\mathcal F}_0$, 
the space
of global sections $H^0
({\mathcal F}\otimes L^{\otimes k})$
has a basis that is defined over $\QQ$. The essence of 
the Conjecture \ref{main} is that the $h$-twisted Eisenstein
series $E_{k,i;h}$ are rational linear combinations 
of the basis elements.  Since $q^0$ coefficients and the $\Gamma_0(l)/\Gamma_1(l)$
transformation properties of $E_{k,i;h}$ are independent 
of the twisting, if the linear span of $E_{k,i;h}$ is defined 
over rationals, then so are the individual elements,
at least in the case of prime $l$. 
Indeed, one can show that in this case the matrix of 
$q^0$ coefficients of $\Gamma_0(l)/\Gamma_1(l)$-translates
of $E_{k,i}$ is nondegenerate and rational. Consequently,
if $Span(E_{k,i;h})$ has a rational basis, $E_{k,i;h}$ 
have rational coefficients in it.

Unfortunately, it is not at all obvious that $Span(E_{k,i;h})$
is defined over $\QQ$. The space of all $h$-twisted
forms of weight $k$ splits up as a direct sum of the space
of cusp forms and its orthogonal complement
under the Petersson pairing. One can hope that 
this splitting respects the rational structures.
In the case 
of prime $l$, the series $E_{k,i;h}(\tau)$ roughly account
for half the dimension of this orthogonal complement.
However, if one could define $E_{k,i;h}(\tau)$ for 
weight $k=1$, then one would expect that all weight one
$h$-twisted modular forms lie in $Span(E_{k,i;h})$,
by an application of  Riemann-Roch formula.

On the plus side, it is quite easy to give
a sufficient condition that assures that the 
Fourier coefficients of $\frac 1{(2\pi\ii)^k}E_{k,i;h}$ are real.
Recall that an $\RR$-structure on the modular curve is given 
by the involution that sends $h(\tau)$ to $\overline {h(-\bar \tau)}$.
Let us calculate $\overline{ E_{k,i;h}(-\bar \tau)}$.
\begin{proposition}
Let $g=\overline{h(-\tau)}$.
For $k\geq 3$ one has
$$
\overline{ E_{k,i;h}(-\bar \tau)}=(-1)^k E_{k,i;g}(\tau)
$$
\end{proposition}

\begin{proof}
and let $G$ be the antiderivative
$\int_{\ii\infty}^\tau g(s)\,ds$ of $g$.
Observe that $\overline{H(-\bar\tau)}=-G(\tau)$. We have
$$
\overline{E_{k,i;h}(-\bar \tau)}=
\sum_{c\in l\ZZ}\sum_{d\in\ZZ,{\rm gcd}(l,d)=1}
(-c\tau+d)^{-k} \ee^{2\pi\ii\frac {di}l} \ee^{-2\pi\ii \Re H(-\frac dc)}
$$
$$
=
\sum_{c\in l\ZZ}\sum_{d\in\ZZ,{\rm gcd}(l,d)=1}
(-c\tau+d)^{-k} \ee^{2\pi\ii\frac {di}l} \ee^{2\pi\ii \Re G(\frac dc)}
=E_{k,-i;g}(\tau)
$$
where we have switched from $c$ to $-c$ in the summation.
It remains to observe $E_{k,-i;g}=(-1)^k E_{k,i;g}$, in view
of the change from $(c,d)$ to $(-c,-d)$ in the summation.
\end{proof}

\begin{corollary}
If all Fourier coefficients of $h(\tau)$ are real, then 
all Fourier coefficients of $\frac 1{(2\pi\ii)^k}E_{k,i;h}(\tau)$ are real.
\end{corollary}

\begin{proof}
Notice that If 
$F(\tau)=\sum_m r_mq^m$ then $\overline{F(-\bar \tau)}=\sum_m
\overline{r_m} q^m$. Then use the above proposition with $h=g$. 
\end{proof}

 \section{Open problems}\label{sec6}
There are numerous open problems, whose solution
would greatly enhance the understanding of the 
$h$-twisted modular forms. We list several of them.

\begin{itemize}

\item
Can one efficiently calculate Fourier
coefficients of an $h$-twisted
Eisenstein series? The problem is that the definition
only provides a relatively slowly convergent series.
Lack of explicit coefficients makes it impossible
to do any kind of computer experiments to uncover
relations among the series or to test Conjecture \ref{main}.

\item
Can one define $h$-twisted Eisenstein series for 
weights $k=1$ or $2$? What will happen to their
modularity properties? The similar question in
the untwisted case is answered by looking at the 
Fourier expansions, but this approach seems to fail
in general.

\item
In the untwisted case, these Eisenstein series
can be constructed from logarithmic derivatives of 
the theta function, see \cite{BG1}. Is there a twisted
analog of the theta function, with the same properties?
In a related comment,
the notion of $h$-twisting can be extended
to the Jacobi forms.

\item
Do any twisted forms have nice infinite product expansions?

\item
What is the meaning of Hecke operators in this setting?
It is easy to see that the usual definition of the Hecke
operators does not produce operators on the spaces
of $h$-twisted forms. Also, if $h$ is a Hecke eigenform,
does this lead to any properties of $h$-twisted forms?
A related question is what should be the analog of 
Dirichlet series and Euler expansions in the $h$-twisted case.

\item
What is the meaning of the Fricke involution in this setting?

\item
What would be an analog of Rankin-Selberg method 
for the Petersson inner product of a usual cusp form 
and product of two twisted Eisenstein series of opposite
twists? The argument of Proposition \ref{orthogonal} shows that
the unfolding trick generally works, but the lack of
Euler expansions is a major hindrance.

\item 
Is it possible to see what the twisted Eisenstein series are
in the case of the divisors coming from the 
Heegner points?

\item
What is the relation between the Eisenstein series of 
opposite twists?

\item
In the untwisted case, there are quadratic relations on
the Eisenstein series that mimic the relations on modular 
symbols. Is there any analog of such relations in the 
twisted case? For instance, dimension counts show
that there are a lot of linear relations between products
of $h$-twisted and untwisted Eisenstein series. Can
one write at least some of them explicitly?
Together with the rationality conjecture,
this might provide an approach to the Birch-Swinnerton-Dyer
conjecture, see \cite{Tate}, 
although at this point in time this is, at best, a very long shot.

\end{itemize}

\end{document}